\documentclass[11pt]{amsart}

\usepackage{amsmath}
\usepackage{amssymb}
\usepackage{graphicx}

\newtheorem{theorem}{Theorem}[section]

\theoremstyle{definition}
\newtheorem{definition}[theorem]{Definition}

\newtheorem{rem}[theorem]{Remark}
\newtheorem{thm}{Theorem}[section]
\newtheorem{lem}[theorem]{Lemma}
\newtheorem{prop}[theorem]{Proposition}

\numberwithin{equation}{section}

\DeclareMathOperator{\diam}{diam}

\newcommand{\N}{\mathbb{N}}

\newcommand{\R}{\mathbb{R}}

\def\rr{{\mathbb R}}

\def\loc{{\mathop\mathrm{\,loc\,}}}

\def\bint{{\ifinner\rlap{\bf\kern.35em--}
		\int\else\rlap{\bf\kern.45em--}\int\fi}\ignorespaces}

\def\bbint{{\ifinner\rlap{\bf\kern.35em--}
		\hspace{0.078cm}\int\else\rlap{\bf\kern.45em--}\int\fi}\ignorespaces}

\newcommand{\slmap}[2]{W^{#1,#2}_\text{loc}(\ball,\R^n)}

\newcommand{\ball}{B^n(0,1)}

\newcommand{\abs}[1]{\vert{#1}\vert}

\newcommand{\dia}[1]{\text{diam}(#1)}

\def\diam{{\mathop\mathrm{\,diam\,}}}

\title[A note to ``Radial limits of quasiregular Local Homeomorphisms'']
{A note to ``Radial limits of quasiregular Local Homeomorphisms''}

\dedicatory{Dedicated to the memory of Yurii Reshetnyak}

\author[C.-Y. Guo]{Chang-Yu Guo}

\address[C.-Y. Guo]{Research Center for Mathematics, Shandong University 266237,  Qingdao, P. R. China and  Frontiers Science Center for Nonlinear Expectations, Ministry of Education, P. R. China}
\email{{\tt changyu.guo@sdu.edu.cn}}

\author[Y. Xuan]{Yi Xuan}
\address[Y. Xuan]{HLM, Academy of Mathematics and Systems Science,
	Chinese Academy of Sciences, Beijing, 100190, P. R. China and School of
	Mathematical Sciences, University of Chinese Academy of Sciences,
	Beijing, 100049, P. R. China}
\email{\tt xuanyi@amss.ac.cn}



\thanks{C.-Y. Guo is supported by the Young Scientist Program of the Ministry of Science and Technology of China (No.~2021YFA1002200), the National Natural Science Foundation of China (No.~12101362) and the Natural Science Foundation of Shandong Province (No.~ZR2022YQ01, ZR2021QA003). The corresponding author Y. Xuan is supported by Key Research of Frontier Sciences, CAS (Grant No.~ZDBS-LY-7002).}

\begin{document}

\begin{abstract}
In this short note, we consider quasiregular local homeomorphisms on uniform domains. We prove that such mappings always can be extended to some boundary points along John curves, which extends the corresponding result of Rajala [Amer. J. Math. 2008].
\end{abstract}

\maketitle


\section{Introduction}
Let $\Omega$ be a domain in $\rr^n$ ($n\geq 2$). A mapping $f\colon \Omega\rightarrow \rr^n$ is called  $K$-quasiregular if it belongs to $W^{1,n}_{\loc}(\Omega,\rr^n)$ and satisfies
\[\|Df(x)\|^n\leq KJ_f(x),\]
for almost every $x\in \Omega$, where $Df(x)$ is the differential of $f$ at $x$, $\|Df(x)\|$ is the operator norm of $Df(x)$ and $J_f(x)$ is the Jacabian determinant of $Df(x)$. If $f$ is additionally assumed to be injective, then we recover the class of $K$-quasiconformal mappings. The theory of quasiconformal mappings in higher dimensions was initiated in earnest by Reshetnyak, Gehring and V\"ais\"al\"a (see \cite{M14} for a nice exposition on the history) and it was successfully applied by Mostow \cite{M68} in the proof of his celebrated rigidity theorem. The further extension to the noninjective quasiregular case was again initiated by Reshetnyak \cite{Re67,Re68} and the basic theory was comprehensively developed by Martio, Rickman and V\"ais\"al\"a in their seminal works \cite{mrv69,mrv70,mrv71}. One of the most profound results for the higher dimensional theory of quasiregular mappings is the following theorem of Reshetnyak \cite{Re67,Re68}: \emph{nonconstant quasiregular mappings are continuous, discrete, open and sense-preserving}. For more on the higher dimensional theory of quasiconformal mappings and quasiregular mappings, see the excellent monographs \cite{Re89,V71,R93}.

In this short note, we shall consider the boundary behaviour of quasiregular mappings. Recall that a classical theorem of Fatou \cite{F06} states that, for a bounded analytic function defined on the unit disk, the radial limits exist at almost every points on the unit circle. It remains, however, open up to now whether one can extend this result to the case of quasiregular mappings in higher dimensions. In fact, it is unknown even for bounded quasiregular mappings. A major breakthrough on this problem was made by Kai Rajala \cite{R08}, who proved that a quasiregular local homeomorphism $f\colon B^n(0,1)\to \R^n$, $n\geq 3$, has radial limits at certain point on the boundary $S^{n-1}(0,1)=\partial B^n(0,1)$. More precisely, there exists a point $\xi\in S^{n-1}(0,1)$ so that the radial limit
\[
\lim_{t\to 1}f(t\xi)=b_{\xi}\in \R^n\qquad \text{exists}.
\]
If we add more restrictions on the growth of the multiplicity function of $f$, then one can prove the almost everywhere existence of radial limits on $S^{n-1}(0,1)$ and even with a Hausdorff dimensional estimate on the ``bad" set, where the radial limits do not exist; see \cite{MR72} and \cite{KMV96}.

We shall restrict our attension to Rajala's result \cite{R08} and consider a natural extension by relaxing the definition domain $B^{n}(0,1)$ to the more general class of bounded uniform domains introduced by Gehring and Osgood \cite{GO79}.

Recall that a bounded domain $\Omega\subset\rr^n$ is called a $c$-uniform domain if each pair of points $x_1,x_2\in \Omega$ can be joined by a curve $\gamma\subset \Omega$ with the following two properties:
\[\ell(\gamma)\leq c|x_1-x_2|,\]
\[\min_{i=1,2}\ell(\gamma(x_i,x))\leq cd(x,\partial\Omega) \text{ for all } x\in \gamma.\]
We shall indeed use a variant of uniform domain introduced by \"Akkinen and Guo \cite{AG2017}. For it, we need also the definition of John domains. Recall that
a bounded domain $\Omega\subset \R^n$ is called a $c$-John domain with center $x_0\in\Omega$ if every point $x\in\Omega$ can be joined by a rectifiable curve $\gamma\colon [0,l]\rightarrow \Omega$, parameterized by arc length, such that $\gamma(0)=x$ and $\gamma(l)=x_0$ with the property that
\begin{equation}\label{eq:John domain}
	d(\gamma(t),\partial\Omega)\geq \frac{t}{c} \text{ for every }t\in [0,l].
\end{equation}
It is easy to check that each boundary point of a John domain $\Omega\subset \R^n$ can be connected to $x_0$ by a curve $\gamma$ satisfying \eqref{eq:John domain}. A curve $\gamma$ with such property is called a $c$-John curve. For a point $\xi\in\partial \Omega$, we shall denote by $I_c(\xi,x_0)$ the collection of all $c$-John curves connecting $\xi$ to $x_0$. Notice that for a bounded $c$-uniform domain $\Omega\in\rr^n$, there exists a point $x_0\in\Omega$ such that $\Omega$ is a $c_1$-John domain with center $x_0$, where $c_1$ depending only on $c$; see \cite[Section 2.17]{V88}.
Thus, we define
\begin{definition}[Uniform domain with center]\label{def:uniform domain with center}
	We say that a bounded domain $\Omega\subset \R^n$ is a $c$-uniform domain with center $x_0$ if it is a uniform domain and at the same time it is a $c$-John domain with center $x_0$.
\end{definition}
Boundary behaviours of quasiregular mappings or even the more general class of mappings of finite distoriton, with certain regularity restrictions, were studied in uniform domains with centers by \"Akkinen and Guo in \cite{AG2017}. In particular, existence and uniqueness of limits along John curves were obtained for these classes of mappings.

From now on, unless specified, the dimension $n$ is always assumed to be bigger than or equal to $3$. Our main result of this short note is the following theorem, which can be regarded as a natural extension of \cite[Theorem 1.1]{R08}.
\begin{thm}\label{thm:main result}
	Let $\Omega\subset\rr^n$ be a $c_0$-uniform domain with center $x_0$ and $f\colon \Omega\rightarrow \rr^n$ a quasiregular local homeomorphism. Then, there exists a point $\xi \in\partial\Omega$ such that for any $c_0$-John curve $\gamma\in I_{c_0}(\xi,x_0)$,
	\[
	\lim_{t\to 0+}f(\gamma(t))=b_{\xi}\in \R^n
	\]
	exists and the limit $b_{\xi}$ is independent of the choice of the $c_0$-John curve $\gamma\in I_{c_0}(\xi,x_0)$.
\end{thm}
As was pointed out earlier, if we add more regularity restrictions on $f$, then we could prove the stronger almost every existence (even possibly with a dimensional estimate) of limits along John curves; see \cite[Section 3]{AG2017} for more along this direction.
Similar to \cite{R08}, we will prove a more general result; see Theorem~\ref{3.1} below. In particular, it would imply that limits along John curves exist for infinitely many points on $\partial\Omega$; see the remark after Theorem~\ref{3.1}.

For spatial $K$-quasiregular mappings, the local homeomorphism assumption in Theorem \ref{thm:main result} can be replaced by a smallness assumption $K\leq 1+\epsilon(n)$; see \cite{mrv71,R05}. The assumption $n\geq 3$ is necessary for the conclusion of Theorem \ref{thm:main result}; see \cite[Page 271]{R08} for a counter-example in the planar case.

The idea for the proof of Theorem \ref{thm:main result} is very similar to the one used by Rajala \cite{R08} module some minor technical modifications: we start from a line segment $\gamma$ on the image, such that $\gamma$ terminates at a point $y$ belonging to the cluster set of $\partial\Omega$, and it has a lift $\gamma^{\prime}$ which approaches a point $x$ on $\partial\Omega$. Choose a John curve $\hat{\gamma}$ connecting $x_0$ to $x$. Then we try to relate the properties of $\gamma'$ to the behaviour of the image of the John curve $\hat{\gamma}$. An essential difference from the situation of \cite{R08} is that the lift $\gamma'$ does not have a simple tangential/non-tangential behaviour. Instead, we make use of suitable quasihyperbolic geodesics and do continuity estimate with the aid of quasihyperbolic distances. Fortunately, the necessary auxiliary results, proved in \cite{R08}, remain valid in our situation so that we may estimate the conformal modulus of certain curve families as in \cite{R08}.

The structure of this paper is as follows. In Section \ref{sec:preliminaries}, we give some preliminaries and fix the notation used in this paper. In Section \ref{sec:basic setting}, we fix the basic setting following \cite{R08} and formulate the more general Theorem~\ref{3.1}. We also collect the main auxillary results from \cite{R08} that are needed for the proof in the next section. As the arguments are almost identical to that used in \cite{R08}, we only briefly outline it and point out the part that needs modification. The proof of Theorem \ref{3.1} is given in Section \ref{sec:proof of main result}. In the final section, we consider the existence of radial limits for the more general class of mappings of $p$-integrable distortion, extending the corresponding result of \cite{A14}.
\medskip

\textbf{Acknowledgement:} We would like to thank Prof.~Jinsong Liu and Prof.~Tuomo \"Akkinen for their interests in this work and for many useful discussions during the preparation of this work.

\section{Preliminaries}\label{sec:preliminaries}
We denote by $B^n(x,r)$ or simply $B(x,r)$ an Euclidean ball in $\rr^n$ with center $x$ and radius $r$ and write $S^{n-1}(x,r)=\partial B^n(x,r)$. The Euclidean distance between two points $x,y\in\rr^n$ is denoted by $|x-y|$, while the Euclidean distance between two subsets $A,B\subset\rr^n$ is $d(A,B)$. For a point $a\in\rr^n$, $d(\{a\},B)$ is abbreviated as $d(a,B)$. The Euclidean diameter of a set $A\subset\rr^n$ is denoted by $\diam(A)$. For a set $A\subset\rr^n$ and an open set $G\subset \rr^n$, we write $A\subset\subset G$, if there exists a compact set $K$ such that $A\subset K\subset G$. Define the spherical cap $C(z, \alpha,w)$ by
\[C(z,\alpha,w)=\{x\in\rr^n: |x-z|=|w-z|,(w-z)\cdot(x-z) >{|x-z|}^2\cos\alpha\}.\]
Moreover, for any $z\in\rr^n$, we define the cone $K(z,\varphi)$ with vertex $z$ and opening angle $\varphi$ by
\[K(z,\varphi)=\{x\in\rr^n: z\cdot (z-x)>|z||z-x|\cos \varphi\}.\]

A curve $\gamma$ in $\rr^n$ is a continuous mapping from an interval $I$ to $\rr^n$ and the subsurve of $\gamma$ between $x,y\in\gamma$ is denoted by $\gamma(x,y)$. The image of $\gamma$ is denoted by $|\gamma|$ and the length of a curve $\gamma$ is $\ell(\gamma)$.

Let $\gamma\colon [a,b)\rightarrow \rr^n$ be a curve and $f$ be a continuous, discrete and open mapping from a domain $\Omega\subset \rr^n$ to $\rr^n$. Suppose that $x\in f^{-1}(\gamma(a))$.
The maximum $f$-lifting of $\gamma$ starting at $x$ is defined to be a curve $\beta\colon [a,c)\rightarrow \Omega$ such that
\begin{enumerate}
	\item $\beta(a)=x$,
	\item $f\circ\beta=\gamma|_{[a,c)}$,
	\item If $c<c^{\prime}\leq b$, then there does not exist a curve $\alpha:[a,c^{\prime})\rightarrow \Omega$ such that $\beta=\alpha|_{[a,c)}$ and $f\circ\alpha=\gamma|_{[a,c^{\prime})}$.
\end{enumerate}
For a continuous, discrete and open mapping $f$, such a maximum lifting always exists; see \cite[Chapter II Section 3]{R93}.

The quasihyperbolic distance $d_{qh}$ between two points $x,y$ in a proper domain $\Omega \varsubsetneq\R^n$ is
\[d_{qh}(x,y) = \inf\int_{\gamma}\frac{ds}{d(z,\partial\Omega)},\]
where the infimum is with respect to all rectifiable curves $\gamma$ joining $x$ and $y$ in $\Omega$. This distance was first introduced by Gehring and Palka \cite{GP76}. A curve $\gamma$ joining $x$ and $y$ in $\Omega$ which attains the infimum, that is,
\[d_{qh}(x, y)=\int_{\gamma}\frac{ds}{d(z,\partial\Omega)}\]
is called a quasihyperbolic geodesic joining $x$ and $y$. According to \cite[Lemma 1]{GO79}, a quasihyperbolic geodesic connecting any two points of a proper domain $\Omega$ always exists.

Let $\Omega\subset\rr^n$ be a bounded domain. The Whitney decomposition $\mathcal{W}= \mathcal{W}(\Omega)$ is a family of open cubes $Q\subset\Omega$, which are pairwise disjoint and satisfy the following two properties:
\begin{itemize}
	\item $\Omega=\bigcup_{Q\in\mathcal {W}}\overline{Q}$;
	\item $\diam(Q)\leq d(Q,\partial \Omega)\leq 4\diam(Q).$
\end{itemize}
It is well-known that, for any proper domain $\Omega\varsubsetneq\rr^n$, such decomposition always exists.

Next, we introduce the concept of the modulus of a curve family.
Let $\Gamma$ be a family of curves in $\Omega$. The conformal $n$-modulus $M\Gamma$ is defined by
\[M\Gamma=\inf_{\rho\in X(\Gamma)}\int_{\rr^n}{\rho(x)}^ndx,\]
where $X(\Gamma)$ is the set of all Borel functions $\rho\colon \rr^n\rightarrow [0,\infty]$ such that
\[\int_\gamma\rho ds\geq 1\text{ for all }\gamma\in \Gamma.\]
For a $K$-quasiregular mapping $f\colon \Omega\rightarrow \rr^n$ and a family $\Gamma$ of curves in $\Omega$, the following important Poletsky's inequality holds (see \cite[Chapter II Theorem 8.1]{R93}):
\[Mf(\Gamma)\leq K^{n-1}M\Gamma.\]
Following \cite{V71}, we introduce the concept of modulus to families of curves which lie on submanifolds of $\rr^n$. Let $S$ be an $(n-1)$-dimensional smooth submanifold of $\rr^n$ and $\Gamma$ be a curve family on $S$. The $n$-modulus of $\Gamma$ with respect to $S$ is defined as
\[M^S_n(\Gamma)=\inf_{\rho\in X(\Gamma)}\int_S{\rho}^ndm_{n-1},\]
where $X(\Gamma)$ is the set of all Borel functions $\rho\colon S\rightarrow [0,\infty]$ so that $\int_\gamma \rho ds\geq 1$ for all $\gamma\in\Gamma$.

Next, we prove a simple lemma about the upper bound of the modulus of the curve family related to the Whitney decomposition.
\begin{lem}\label{WC}
	Let $\Omega$ be a domain and $\mathcal{W}(\Omega)$ be the Whitney decomposition of $\Omega$. For any cube $Q\in\mathcal{W}(\Omega)$, let $\Gamma_Q$ be the family of all curves in $\Omega$ that connect $\partial\Omega$ and $Q$. Then, there exists a constant $C$ independent of $Q$ such that $M\Gamma_Q\leq C.$
\end{lem}
\begin{proof}
	Fix a cube $Q\in\mathcal{W}(\Omega)$. Denote the center of $Q$ by $x_Q$. Then, by the property of the Whitney decomposition, we have
	\[Q\subset B\big(x_Q,\diam(Q)/2\big),\]
	and
	\[B\big(x,\diam(Q)\big)\subset \Omega.\]
	Thus, by an elementary property of modulus, we have that $M\Gamma_Q\leq M\Gamma_d$, where $\Gamma_d$ is the family of all curves joining $B(x,d/2)$ and $\rr^n\backslash B(x,d)$. It is well-known (see e.g.~\cite{V71}) that \[M\Gamma_d\leq C,\]
	from which we infer that $M\Gamma_Q\leq C$.
\end{proof}

\section{Basic setting and auxillary results}\label{sec:basic setting}
Throughout this section, we assume that $\Omega\subset\rr^n$ is a $c$-uniform domain with center $x_0$ and $f\colon \Omega\to \R^n$ is a quasiregular local homeomorphism. For a $c$-uniform domain $\Omega$ with center $x_0$, the $c$-John curve $\gamma\colon [0,l]\rightarrow \Omega$ connecting $x_0$ to a boundary point $x_1$ is parameterized by arc length with $\gamma(0)=x_1$ and $\gamma(l)=x_0$.  Without loss of generality, we assume that $f(x_0)=0$. Then, by the property of the local homeomorphism, we obtain a small enough $\delta >0$ and a neighborhood $U$ of $x_0$ so that the restriction of
$f$ on $U$ is a homeomorphism from $U$ to $B(0,\delta)$. For any $y\in S^{n-1}(0,1)$, define \[\gamma_y\colon [0,\infty)\rightarrow \rr^n, \gamma_y(t)=ty.\]
Suppose $\tilde{\gamma}_y$ is the maximum $f$-lifting of $\gamma_y$ starting at $x_0$.
For each $y\in S^{n-1}(0,1)$, there exist a number $\lambda(y)\in [\delta,\infty]$, a point $x_y\in \partial\Omega$ and an increasing sequence $\{t_i\}_{i\geq 1}$ satisfying \[t_i\rightarrow\lambda(y) \text{ as } i\rightarrow \infty\] and
\[\lim_{i\rightarrow\infty}\tilde{\gamma}_y(t_i)=x_y.\]
Set a star-like domain $G\subset \R^n$ as
\[G=\{ty:y\in S^{n-1}(0,1),t\in[0,\lambda(y))\}.\]
Denote the $x_0$-component of $f^{-1}(G)$ by $G^{\prime}$ and denote the restriction of $f$ on $G^{\prime}$ by $g$. Then, it is easy to check that $g$ is a homeomorphism (see \cite[Page 275]{R08}). Thus, $g$ is a quasiconformal mapping.
Moreover, set
\[F=\{z_y:y\in S^{n-1}(0,1)\}\backslash\{\infty\},\]
where
\[z_{y}=\left\{\begin{array}{ll}
	\gamma_{y}(\lambda(y)), & \delta \leq \lambda(y)<\infty, \\
	\infty, & \lambda(y)=\infty.
\end{array}\right.\]
Then $F\neq \emptyset$, since otherwise we would get a quasiconformal mapping $g^{-1}\colon \R^n\to \Omega$, which is impossible as quasiconformal image of $\R^n$ is $\R^n$.

We denote
by $H$ the set of all $z_y\in F$ such that there exist two positive constants $r_y$ and $\varphi_y$ satisfying \[K(z_y,\varphi_y)\cap B(z_y,r_y)\subset G.\]
Then by the local homeomorphism property of $f$, we may find a point $z_y\in H$ minimizing $|z|$ over all $z\in H$. Then $z_y\in H$ so that $H\neq \emptyset$.

Now, we are ready to state our main theorem, from which Theorem \ref{thm:main result} follows.
\begin{thm}\label{3.1}
	Suppose $\Omega$ and $f$ are given as above. Then, for any point $y\in S^{n-1}(0,1)$ with $z_y\in H$ and any $c$-John curve $\gamma$ joining $x_0$ and $x_y$ in $\Omega$, we have
	\begin{equation}\label{2}
		\lim_{t\rightarrow 0+}f(\gamma(t))=z_y.
	\end{equation}
\end{thm}

Reasoning as in \cite[Paragraph after Theorem 3.1]{R08}, we can prove that $H$ contains infinitely many elements. Thus, Theorem \ref{3.1} implies that the conclusion of Theorem \ref{thm:main result} holds at infinitely many points on $\partial\Omega$. However, it is not yet clear, even in the setting when $\Omega=B^n(0,1)$, whether one can deduce measure estimates for the size of the set at which the limits exist.

Next, we collect some useful technical lemmas from \cite[Section 3]{R08} that shall be needed in the proof of Theorem \ref{3.1}.
The first one is a topological result, which will be used later to estimate the conformal modulus of certain curve families. For a spherical cap $C\subset S^{n-1}(a,r)$, $\partial C$ represents the relative boundary of $C$ with respect to $S^{n-1}(a,r)$.
\begin{lem}\label{3.2}
	Assume that $z_y\in F$ and that $\hat{\gamma}\colon [0,1]\rightarrow\Omega$ is a curve with $\hat{\gamma}(0)\in|\tilde{\gamma}_y|$. Furthermore, suppose that, for some $r>0$,
	$f(\hat{\gamma}(0))\in B(z_y,r)$ and that the set \[S:=\{s:f(\hat{\gamma}(s))\in S^{n-1}(z_y,r)\}\neq \emptyset.\]
	Set $s^r=\min_{s\in S}.$ Then there exist an angle $0<\alpha\leq\pi$ and a point \[k_r\in\partial C(z_y,\alpha,f(\hat{\gamma}(s^r))),\] such that each curve $\eta$ joining $f(\hat{\gamma}(s^r))$ to $k_r$ in $C(z_y,\alpha,f(\hat{\gamma}(s^r)))$ satisfies that the maximal lift $\eta^{\prime}$ of $\eta$ starting at $\hat{\gamma}(s^r)$ has the property that
	\[\overline{|\eta^{\prime}|}\cap\partial\Omega\neq \emptyset.\]
\end{lem}
\begin{proof}
	The proof is essentially contained in \cite[Proof of Lemma 3.2]{R08}. The argument is basically topological and uses only the local homeomorphism property of $f$. For the convenience of readers, we give a brief outline here.
	
	For $\varphi>0$, denote by $C_{\varphi}^{\prime}$ the component with $\hat{\gamma}\left(s^{r}\right)$ of $f^{-1}\left(C\left(z_{y}, \varphi, f\left(\hat{\gamma}\left(s^{r}\right)\right)\right)\right)$.
	Let $\alpha \in(0, \pi]$ be the maximal angle such that
	$$
	f|_{ C_{\alpha}^{\prime}}\colon C_{\alpha}^{\prime} \rightarrow C\left(z_{y}, \alpha, f\left(\gamma\left(s^{r}\right)\right)\right)
	$$
	is a homeomorphism and let $h=f|_{ C_{\alpha}^{\prime}}$. We consider two cases.
	
	\textbf{Case 1}: $\alpha<\pi$.
	
	We argue by contradiction. If the lemma does not hold, then for each \[p \in \partial C\left(z_{y}, \alpha, f\left(\hat{\gamma}\left(s^{r}\right)\right)\right),\] there is a curve $\eta$ joining $f\left(\hat{\gamma}\left(s^{r}\right)\right)$ and $p$ in $C\left(z_{y}, \alpha, f\left(\hat{\gamma}\left(s^{r}\right)\right)\right)$ with the property that the lift $\eta^{\prime}$ of $\eta$ starting at $\hat{\gamma}\left(s^{r}\right)$ satisfying that \[\left|\eta^{\prime}\right| \subset \Omega_\epsilon\] for some $\epsilon>0,$ where \[\Omega_{\epsilon}=\{x\in\Omega:d(x,\partial\Omega)>\epsilon\}.\]
	By \cite[Chapter III Lemma 3.3]{R93}, $h^{-1}$ extends to a map on a neighborhood of $p$. As $p \in \partial C\left(z_{y}, \alpha, f\left(\hat{\gamma}\left(s^{r}\right)\right)\right)$ is arbitrary, $h^{-1}$ extends to a homeomorphism of $\overline{C}\left(z_{y}, \alpha, f\left(\hat{\gamma}\left(s^{r}\right)\right)\right)$ onto $\overline{C_{\alpha}^{\prime}}$. Moreover, \cite[Chapter III Lemma 3.2]{R93} ensures that we may extend $h^{-1}$ to be a homeomorphism of a neighborhood of $\bar{C}\left(z_{y}, \alpha, f\left(\hat{\gamma}\left(s^{r}\right)\right)\right)$ onto its image. Thus, $\alpha$ is not maximum, which is a contradiction.
	
	\textbf{Case 2}: $\alpha=\pi$.
	
	We denote by $p$ the unique point of $\partial C\left(z_{y}, \pi, f\left(\hat{\gamma}\left(s^{r}\right)\right)\right)$ and argue by contradiction. If the lemma fails, then there is a curve $\eta$ joining $f\left(\hat{\gamma}\left(s^{r}\right)\right)$ and $p$ with the property that the maximal lift $\eta^{\prime}$ of $\eta$ starting at $\hat{\gamma}\left(s^{r}\right)$ satisfying that
	\begin{equation}\label{3.4}
		\left|\eta^{\prime}\right| \subset \Omega_\epsilon
	\end{equation}
	for some $\epsilon>0$. Next, we define a nested sequence of sets as
	$$
	W_{i}=h^{-1}\left(B\left(p, i^{-1}\right) \cap C\left(z_{y}, \pi, f\left(\hat{\gamma}\left(s^{r}\right)\right)\right)\right)
	$$
	As $n \geq 3$, each $W_{i}$ is connected and $\overline{W_{i}}$ is also connected. Thus,
	$$
	\mathcal{K}=\bigcap_{i=1}^{\infty} \overline{W_{i}}
	$$
	is compact and connected. Thus, by \eqref{3.4}, we obtain that
	$$
	\mathcal{K} \cap \Omega \neq \emptyset
	$$
	However, as $f$ is discrete and $$
	\mathcal{K} \cap \Omega \subset f^{-1}(p),
	$$
	$\mathcal{K}$ has to be a singleton.
	
	The remaining argument is exactly the same as in the proof of Lemma 3.2 of \cite{R08}. One uses a topological argument to construct a curve $\beta$ that has two lifts, both starting at $\hat{\gamma}(s^r)$, contradicting with the local homeomorphism property of $f$.
\end{proof}

In the following, we fix a point $z_y\in H$ and set $\varphi_y=\varphi_0$ and $r_y=r_0$. Without loss of generality, assume that $r_0\leq \frac{1}{2}$ and $|x_y-x_0|=1$.
For $0<r<1$, set
\[t_r=\min\{t\in[0,\lambda(y)):\tilde{\gamma}_y(t)\in S^{n-1}(x_y,r)\}\]
and
\begin{equation}\label{101}
	x_r=\tilde{\gamma}_y(t_r)\in S^{n-1}(x_y,r).
\end{equation}
Then, it is clear that
\[
|g(x_{r_1})-z_y|<|g(x_{r_2})-z_y|,\text{ whenever }r_1<r_2.
\]
Following \cite{R08}, we define a function $\theta\colon (0,1]\rightarrow (0,1]$ as
$$\theta(t)=d(x_t,\partial\Omega).$$
Then, it is easy to check that $\theta$ is right continuous. Furthermore, we have the following integral estimate on $\theta$.
\begin{lem}[Proposition 4.1, \cite{R08}]\label{3.3}
	If $f(x_s)\in B(z_y,r_0/2)$, then, for any $0<t<s/2$, we have
	\[\int_t^s\frac{dr}{\theta(r)}\leq C_0\log\frac{C_1}{|f(x_t)-z_y|},\]
	where two positive constants $C_0,C_1$ only depend on $\varphi_0$, $r_0$, $K$ and $n$.
\end{lem}
The definition  of $\theta$ as in the above lemma differs from that used in \cite[Proposition 4.1]{R08}. However, the proof there only uses the quantitative equivalence of quasiconformality with local quasisymmetry and thus applies to our setting without modification.


Finally, we recall the following technical function lemma about $\theta$.
\begin{lem}[Lemma 5.1, \cite{R08}]\label{5.1}
	Let $\theta\colon (0, 1]\rightarrow (0, 1]$ be a right continuous function such
	that $t_i\rightarrow 0$ if $\theta(t_i)\rightarrow 0$, and with the property that
	\begin{equation}\label{unsure}\lim_{t\rightarrow 0}\theta (t)t^{-1}=0.\end{equation}
	Then, for any $\epsilon>0$, there is a number $M$ and a decreasing sequence $T_i\rightarrow 0$ so that
	\[\log \frac{1}{\theta (T_i)} \leq \epsilon \int_{T_i}^1 \frac{dr}{\theta(r)}+M\text{ for any } i\geq 1.\]
	
\end{lem}

\section{Proof of the Main Result}\label{sec:proof of main result}

\begin{proof}[Proof of Theorem~\ref{3.1}]
	Following \cite[Proof of Theorem 3.1]{R08}, we prove Theorem \ref{3.1} by a contradiction argument.
	Suppose the conclusion of Theorem \ref{3.1} does not hold true. Then there exist a $c$-John curve $\gamma$ connecting $x_0$ and $x_y$ in $\Omega$, a constant $m>0$ and a decreasing sequence $t_j\rightarrow 0$, such that
	\[|f(\gamma (t_j))-z_y|\geq m \text{ for all }j\geq 1.\]
	Set $a_j=\gamma(t_j)$ and fix a constant $s$ as in Lemma~\ref{3.3}. Since
	\[|x_y-a_j|\leq \ell(\gamma(x_y,a_j))=t_j,\]
	we may assume that
	\[|x_y-a_j|<s\text{ for all } j\geq 1.\]
	For each $j\geq 1$, we choose a point $b_j=x_{\alpha(j)}$, for some
	$0<\alpha(j)<s$, so that
	\[|b_j-x_y|<|a_j-x_y|,\]
	where $x_{\alpha(j)}$ is given as in \eqref{101} and we shall determine $\alpha(j)$ later. Then, we have
	\begin{equation}\label{ww}
		\alpha(j)=|b_j-x_y|\xrightarrow{j\to \infty} 0.
	\end{equation}
	
	Let $\beta$ be a quasihyperbolic geodesic joining $a_j$ and $b_j$ and denote by $\mathcal{W}(\Omega)$ the Whitney decomposition of $\Omega$. The set of all cubes in $\mathcal{W}(\Omega)$ that intersects the quasihyperbolic geodesic $\beta$ is labeled as $\{Q_i\}_{1\leq i\leq N_j}$. Then we have
	\begin{equation}\label{ggg}
		N_j\leq Cd_{qh}(a_i,b_i),
	\end{equation}
	where $C$ is a constant only depending on $n$ (see e.g.~\cite[Page 421]{KOT01}). As $\Omega$ is a $c$-uniform domain, we have (see e.g.~\cite{GO79})
	\begin{equation}\label{hhh}
		d_{qh}(a_j,b_j)\leq C\log\Big(1+\frac{|a_j-b_j|}{\min\{d(a_j,\partial\Omega),d(b_j,\partial\Omega)\}}\Big),
	\end{equation}
	where $C$ is a constant only depending on $n$ and $c$. Since $\gamma$ is a $c$-John curve, we have
	\[d(a_j,\partial\Omega)\geq \ell(\gamma(a_j,x_y))/c\geq |a_j-x_y|/c\geq \theta(\alpha(j))/c=d(b_j,\partial\Omega)/c.\]
	Furthermore,
	\[|a_j-b_j|\leq |a_j-x_y|+|b_j-x_y|\leq 2|x_y-a_j|.\]
	Thus, it follows from \eqref{ggg}, \eqref{hhh} and above estimates that
	\begin{equation}\label{sss}N_j\leq C\log(1+C\frac{|x_y-a_j|}{\theta(\alpha(j))}).
	\end{equation}
	
	Set
	\begin{equation}\label{ooo}
		r_j=|f(b_j)-z_y|\quad \text{and}\quad s_j=|f(a_j)-z_y|\geq m.
	\end{equation}
	
	\textbf{Claim:} $r_j\rightarrow 0$ as $j \rightarrow \infty$.
	
	We prove the claim by contradiction. If the claim fails, then there exists a constant $\epsilon$ and a subsequence $b_{k_j}$ such that $|f(b_{k_j})-z_y|>\epsilon$. As $f(b_{k_j})\in \gamma_y([0,\lambda(y)))$ for all $j\geq 1$, it follows that \[\{f(b_{k_j})\}\subset\subset G.\] Since $g$ is a homeomorphism, we have \[\{b_{k_j}\}\subset \subset G^{\prime},\]
	which is a contradiction to \eqref{ww}. Therefore, the claim is true.
	
	Then, for all $r\in (r_j,s_j)$, we have
	\[f(b_j)\in B(z_y,r),\]
	and
	\[f(a_j)\in \rr^n\backslash B(z_y,r).\]
	For each $r\in (r_j,s_j)$, we may choose a parametrization of $\beta\colon [0,1]\rightarrow \Omega$ with $\beta(0)=b_j$ and $\beta(1)=a_j$ so that the assumptions of Lemma~\ref{3.2} are satisfied.
	Set $p_r=f(\beta(s^r))$, where $s^r$ is given as in Lemma~\ref{3.2}.
	Then, Lemma \ref{3.2} implies that there exist a spherical cap $C_r$ and a point $k_r\in \partial C_r$ such that, if we denote by $\Gamma_r$ the set of all curves connecting $p_r$ and $k_r$ in $C_r$, then the maximum lift $\gamma^{\prime}$ of each $\gamma\in\Gamma_r$ starting at $\beta(s^r)$ satisfies
	\begin{equation}\label{ttt}\overline{|\gamma^{\prime}|}\cap \partial\Omega\neq\emptyset.
	\end{equation}
	An application of \cite[Theorem 10.2]{V71} implies that
	\[M^n_S\Gamma_r\geq \frac{C}{r},\]
	where the positive constant $C$ depending only on $n$. Thus, for the curve family \[\Gamma =\bigcup_{r\in(r_j,s_j)}\{\gamma:\gamma\in \Gamma_r\},\]
	integrating the previous estimate with respect to $r$ gives
	\begin{equation}\label{zzz}
		M\Gamma\geq C\log \frac{s_j}{r_j}.
	\end{equation}
	
	Next, we set
	\[\Gamma^{\prime}=\{\gamma^{\prime}:\gamma^{\prime}\text{ is the maximum lifting of some }\gamma\in\Gamma_r \text{ starting at } \beta(s^r),r\in(r_j,s_j)\}.\]
	Remember that the cubes, in the Whitney decomposition, that intersect $\beta$ are denoted by $\{Q_i:1\leq i\leq N_j\}$. Denote by $\Gamma_i$ the set of all curves connecting $Q_i$ and $\partial\Omega$, where $1\leq i \leq N_j$. Then, by Lemma~\ref{WC},
	\[M\Gamma_i\leq C,\text{ for all }1\leq i\leq N_j.\]
	By \eqref{ttt} and the subadditivity of the modulus,
	it follows that
	\[M\Gamma^{\prime}\leq \sum_{i=1}^{N_j}M\Gamma_i\leq CN_j.\]
	Combining the above estimate with \eqref{sss}, we get
	\[M\Gamma^{\prime}\leq C\log(1+C\frac{|x_y-a_j|}{\theta(\alpha(j))}).\]
	It follows then from Poletsky's inequality, the above inequality, \eqref{zzz} and \eqref{ooo} that
	\begin{equation}\label{cc}\log m\leq \log s_j=\log\frac{s_j}{r_j}-\log\frac{1}{r_j}\leq C_2\log(1+C_3\frac{|x_y-a_j|}{\theta(\alpha(j))})-\log\frac{1}{r_j}.\end{equation}
	
	Notice that the function $\theta$ satisfies all the requirements of Lemma~\ref{5.1} except \eqref{unsure}. Next, we consider two cases depending on whether $\theta$ satisfies \eqref{unsure}.
	\medskip
	
	\textbf{Case 1}: \eqref{unsure} does not hold.
	\medskip
	
	It follows that
	$$
	\limsup _{t \rightarrow 0} \theta(t) t^{-1} \geq \tau^{\prime}>0 .
	$$
	For now, fix $$\tau\in(0,\tau^{\prime}).$$ Set
	$$
	m_{j}=\left|x_{y}-a_{j}\right|<s.
	$$
	Define $\epsilon_{j} m_{j}$ to be the supremum of all $r \in (0,m_{j}]$ satisfying that
	$\theta(r) \geq \tau r.$
	It is easy to check that, $\epsilon_{j} m_{j}>0$ for all $j \geq 1$. Thus, for any $j\geq 1$, there is a number $\alpha(j) \in\left[\epsilon_{j} m_{j} / 2, \epsilon_{j} m_{j}\right]$ such that $\theta(\alpha(j)) \geq \tau \alpha(j)$.
	
	Next, if we assume that $\epsilon_{j}>\epsilon>0$ for all $j \geq 1$, then \eqref{cc} implies that
	$$
	\log m \leq C\log\left( C\frac{2 m_{j}}{\tau \epsilon m_{j}}+1\right)-\log \frac{1}{r_{j}} \leq C(\tau, \epsilon)-\log \frac{1}{r_{j}} \rightarrow-\infty
	$$
	as $j \rightarrow \infty$. This is a contradiction.
	
	Thus, by taking a subsequence if necessary, we may assume that $\epsilon_{j}$  tends to $0$ as $j \rightarrow \infty$ and $\epsilon_{j}$ is decreasing.
	Now, as $\theta(r)<\tau r$ for all $r \in\left[\epsilon_{j} m_{j}, m_{j}\right]$, by Lemma~\ref{3.3}, we obtain
	$$
	\log \frac{C}{r_{j}} \geq C^{-1} \int_{\epsilon_{j} m_{j}}^{m_{j}} \frac{d r}{\theta(r)} \geq \frac{\log \frac{1}{\epsilon_{j}}}{C \tau}.
	$$
	Then, the above inequality and \eqref{cc} yield
	$$
	\log m \leq C\log\left( \frac{2C}{\tau}+\epsilon_j\right)+(C_1-\tau^{-1} C_2^{-1}) \log \frac{1}{\epsilon_{j}}+C.
	$$
	Choose $\tau$ small enough so that the coefficient $C_1-\tau^{-1}C_2^{-1}<0$.  Then, the right hand side of the above inequality tends to $-\infty$ as $j \rightarrow \infty$. This is a contradiction and we finish the proof of this case.
	\medskip
	
	\textbf{Case 2}: \eqref{unsure} holds.
	\medskip
	
	We apply Lemma~\ref{5.1} with $\epsilon=C_0^{-1}C_2^{-1}$. Choose, for each $j\geq 1$, a number $i=i(j)$ such that \[T_{i(j)}=|x_{T_{i(j)}}-x_y|<|a_j-x_y|,\] where $T_{i}$ is as in Lemma~\ref{5.1}. Set $\alpha(j)=T_{i(j)}$. Then, by Lemma~\ref{3.3} and Lemma~\ref{5.1}, we have
	$$
	\begin{aligned}
		C_{2} \log \frac{1}{\theta\left(T_{i}\right)}-\log \frac{C_1}{r_{j}} & \leq C_{2} \log \frac{1}{\theta\left(T_{i}\right)}-C_{0}^{-1} \int_{T_{i}}^{s} \frac{d r}{\theta(r)} \\
		& \leq M C_{2}+C_{0}^{-1} \int_{s}^{1} \frac{d r}{\theta(r)}.
	\end{aligned}
	$$
	Thus, by combining the above inequality and \eqref{cc} we infer that
	$$
	\log m \leq C_{2}\left(\log (C_3\left|x_{y}-a_{j}\right|+\theta(T_i))+M\right)+\log C_1+C_{0}^{-1} \int_{s}^{1} \frac{d r}{\theta(r)} \xrightarrow{j\to\infty}-\infty,
	$$
	which is a contradiction. Therefore, the proof of Theorem~\ref{3.1} is complete.
	
\end{proof}

\section{Further discussion}\label{sec:further discussion}
In this section, we give a brief discussion on the existence of radial limits for certain mappings of finite distortion with extra growth conditions on the Jacobian determinant, along the line of \cite{A14}. For this, we first introduce the definition of mappings of $p$-integrable distortion; for more on the theory of mappings of finite distortion, see the monographs \cite{IM01,HK14}.


Throughout this section, we assume that $f\in\slmap{1}{1}$, $J_f\in L^1_\text{loc}(\ball)$
and that there exists a measurable function $K(x)\geq 1$, finite almost everywhere, such that $f$ satisfies the distortion inequality
\begin{equation}\label{dist}
	\abs{Df(x)}^n\leq K(x) J_f(x)
\end{equation}
for almost every $x\in\ball$. Moreover we assume that
\[\frac{\abs{Df(x)}^n}{\log(e+\abs{Df(x)})}\in L^1_\text{loc}(\ball).\]

\begin{definition}\label{def:bounded p mean distortion}
	We say that $f\colon B^n(0,1)\to \R^n$ is a mapping of  $p$-integrable distortion if there exists $M<\infty$ such that
	\[\int_{B^n(0,1)}K^p(x)\,dx\leq M.\]
\end{definition}

Our main result of this section is the following result that extends \cite[Corollary 1.2]{A14}.
\begin{thm}\label{thm:p-mean distortion}
	For $p>n-1$, let $f\colon\ball\to\R^n$ be a mapping of $p$-integrable distortion and further assume that the Jacobian determinant of $f$ satisfies growth condition
	\[\int_{B(0,r)}J_f(x)\,dx\leq c(1-r)^{-a}\]
	for some $a\in[0,n-1-n/p)$, then
	\[\dim_{\mathcal{H}}{E(f)}\leq a+\frac{n}{p},\]
	where $\dim_{\mathcal{H}}$ refers to the Hausdorff dimension and $E(f)$ consists of all points $x\in S^{n-1}(0,1)$ so that $f$ does not have a radial limit at $x$. In particular, we have radial limits almost everywhere in $\partial\ball$.
\end{thm}

Fix $p>n-1$ and set
\[p^*=\frac{np}{p+1}.\]
Then it is clear that $n-1<p^*<n$ and $p=p^*/(n-p^*)$. We shall need the following standard continuity estimate.
\begin{prop}\label{SI}
	Let $f\colon \ball\rightarrow\R^n$ be a mapping of $p$-integrable distortion with $p>n-1$. If $B(z,r)$ is such that $B(z,2r)\subset\ball$, then for every $x,y\in B(z,r)$
	\[\abs{f(x)-f(y)}^{p^*}\leq C(n,p^*)t^{p^*-n+1}\int_{\partial B(z,t)}\abs{Df(x)}^{p^*}\,dx,\]
	for almost every $t\in(r,2r)$.
\end{prop}
\begin{proof}
	First note that $f\in W_\text{loc}^{1,p^*}(\ball,\R^n)$ as $p^*<n$. Proposition~\ref{SI} follows now directly from the following  well-known Sobolev embedding on spheres (see e.g.~\cite[Lemma 3.2]{A14})
	\[\abs{f(x)-f(y)}^s\leq C(n,s)t^{s-n+1}\int_{\partial B(z,t)}\abs{Df(x)}^s\,dx.\]
	and the fact that continuous and open map in the Euclidean space is monotone, i.e.
	\[\dia{f(B^n(x,r))}\leq \dia{f(S^{n-1}(x,r))}\]
	for all $B^n(x,r)\subset\ball$.
\end{proof}

Now we are ready to prove Theorem \ref{thm:p-mean distortion}
\begin{proof}[Proof of Theorem \ref{thm:p-mean distortion}]
	The proof is similar to \cite[Proof of Theorem 1.1]{A14}.
	Define $\gamma_y\colon [0,1]\rightarrow \ball$ to be the radial segment associated with $y\in S^{n-1}(0,1)$. Then, for each $k=1,2,\ldots$ set
	\[\gamma_{y,k}=\gamma_y|_{[r_k,r_{k+1}]},\]
	where $r_k=1-2^{-k}$. Moreover, set
	\[A_k=\{y\in S^{n-1}(0,1) : \dia{f\abs{\gamma_{y,k}}}\geq k^{-1-\delta/n}\}.\]
	Now
	\[E(f)\subset\bigcap_{j=1}^\infty\bigcup_{k=j}^\infty A_k,\]
	so we set
	\[E_0=\bigcap_{j=1}^\infty\bigcup_{k=j}^\infty A_k.\]
	
	Next we want to have suitable covers for the sets $A_k$ and this can be done with the aid of Proposition \ref{SI}. Let $y\in S^{n-1}(0,1)$ and $k\in\mathbb{N}$, and define $x_0=\frac{r_{k+1}+r_k}{2}y$ and $R_k=\frac{r_{k+1}-r_k}{2}$. Then $B(x_0,2R_k)\subset B(0,r_{k+3})\setminus B(0,r_{k-1})$, and by Corollary \ref{SI},
	\[\abs{f(x)-f(y)}^{p^*}\leq C(n,p^*)t^{{p^*}-n+1}\int_{\partial B(x_0,t)}\abs{Df(x)}^{p^*}\,dx\]
	for almost every $t\in(R_k,2R_k)$, whenever $x,y\in B(x_0,R_k)$. We integrate this over the range $(R_k,2R_k)$ with respect to variable $t$ to obtain
	\begin{align*}
		\abs{f(x)-f(y)}^{p^*} & \leq C(n,p^*)R_k^{p^*-n}\int_{B(x_0,2R_k)}\abs{Df(x)}^{p^*}\,dx.
	\end{align*}
	Moreover, we have by continuity of $f$ that
	\begin{equation}\label{Est}
		\dia{f\abs{\gamma_{y,k}}}^{p^*}\leq C(n,{p^*})R_k^{{p^*}-n}\int_{B(x_0,2R_k)}\abs{Df(x)}^{p^*}\,dx.
	\end{equation}
	Now let $y\in A_k$. Then
	\[1\leq C(n,{p^*})(k^{1+\delta/n})^{p^*} R_k^{{p^*}-n}\int_{r_{k-1}}^{r_{k+3}}\int_{S(y,\alpha(k))}\abs{Df(tx)}^{p^*}t^{n-1}\,dx\,dt,\]
	where $\alpha(k)=\arcsin\big(2\frac{R_k}{\abs{x_0}}\big)$ and $S(z,r)$ denotes spherical cap in $S^{n-1}(0,1)$. Then for each $k$ we cover sets $A_k$ with spherical balls $S(y_k,\alpha(k))$ and use Vitali's covering lemma to get a finite subcollection of pairwise disjoint spherical balls, say $\{S(y_i^k,\alpha(k))\}_{i=1}^{p_k}$ such that
	\[A_k\subset\bigcup_{i=1}^{p_k}S(y_i^k,5\alpha(k)).\]
	Notice that
	\[E_0=\bigcap_{j=1}^\infty\bigcup_{k=j}^\infty A_k\subset\bigcup_{k=s}^\infty A_k\]
	for every $s\in\mathbb{N}$. Let $t\in(a+\frac{n}{p},n-1)$ and $\varepsilon>0$. Furthermore, choose $M\in\mathbb{N}$ such that
	\[\dia{S(y_i^k,\alpha(k))}\leq\varepsilon\]
	for all $k\geq M$. Then
	\begin{equation}\label{mofd1}
	\begin{aligned}
		\mathcal{H}_{\varepsilon}^t(E(f)) & \leq\mathcal{H}_{\varepsilon}^t(\bigcup_{k=M}^\infty A_k) \leq\sum_{k=M}^\infty\mathcal{H}_{\varepsilon}^t(A_k)\\
		& \leq\sum_{k=M}^\infty\sum_{i=1}^{p_k} \dia{(S(y_k^{i},5\alpha(k))}^t\\
		& \leq C(n)\sum_{k=M}^\infty  p_k2^{-kt},
	\end{aligned}
	\end{equation}
	since $\alpha(k)\leq C(n)2^{-k}$.
	
	Now we need an estimate for the amount of spherical balls needed to cover each $A_k$. Denote $x_0^{i,k}=\frac{r_{k+1}+r_k}{2}y_i^k$ and $B_i^k=B(x_0^{i,k},2R_k)$.
	Summing both sides in (\ref{Est}) yields
	\begin{align*}
		p_k&\leq C(n,{p^*})(k^{1+\delta/n})^{p^*}R_k^{{p^*}-n}\sum_{i=1}^{p_k}\int_{B_i^k}\abs{Df(x)}^{p^*}\,dx.
	\end{align*}
	Therefore, since $R_k=C2^{-k}$, the distortion inequality \eqref{dist} and H\"older's inequality give
	\begin{align*}
		p_k & \leq C(n,{p^*})(k^{1+\delta/n})^{p^*}2^{k(n-{p^*})}\sum_{i=1}^{p_k}\int_{B_i^k}(\abs{Df(x)}^n)^{{p^*}/n}\,dx\\
		& \leq C(n,{p^*})(k^{1+\delta/n})^{p^*}2^{k(n-{p^*})}\sum_{i=1}^{p_k}\int_{B_i^k}K(x)^{{p^*}/n}J_f(x)^{{p^*}/n}\,dx\\
		& \leq C(n,{p^*})(k^{1+\delta/n})^{p^*}2^{k(n-{p^*})}\sum_{i=1}^{p_k}\Big(\int_{B_i^k}K(x)^{\frac{{p^*}}{n-{p^*}}}\,dx\Big)^{\frac{n-{p^*}}{n}}\Big(\int_{B_i^k}J_f(x)\,dx\Big)^{\frac{p^*}{n}}\\
		& = C(n,{p^*})(k^{1+\delta/n})^{p^*}2^{k(n-{p^*})}\sum_{i=1}^{p_k}\Big(\int_{B_i^k}K(x)^p\,dx\Big)^{\frac{n-{p^*}}{n}}\Big(\int_{B_i^k}J_f(x)\,dx\Big)^{\frac{p^*}{n}}.
	\end{align*}
	Applying the following H\"older's inequality with $\theta=n/p^*$
	\[
	\sum_{i=1}^{p_k}a_i^{\frac{n-{p^*}}{n}}b_i^{\frac{p^*}{n}}\leq \Big(\sum_{i=1}^{p_k}1\Big)^{1-\frac{1}{\theta}}\left(\sum_{i=1}^{p_k}\Big(a_i^{\frac{n-{p^*}}{n}}b_i^{\frac{p^*}{n}}\Big)^{\theta} \right)^{\theta^{-1}}
	\]
	and then raising to power $\theta$ on both sides, we conclude 
	\[p_k\leq C(n,p^*)k^{n+\delta}2^{k\frac{n(n-{p^*})}{p^*}}\sum_{i=1}^{p_k}\Big(\int_{B_i^k}K(x)^p\,dx\Big)^{\frac{1}{p}}\int_{B_i^k}J_f(x)\,dx.\]
	By $p$-integrability of the distortion we have
	\[\int_{B_i^k}K(x)^p\,dx\leq\int_{B^n(0,r_{k+3})}K(x)^p\,dx\leq M.\]
	Thus by assumption on the Jacobian determinant
	\begin{align*}
		p_k&\leq C(n,p^*,M) k^{n+\delta}2^{k\frac{n(n-{p^*})}{p^*}}\sum_{i=1}^{p_k}\int_{B_i^k}J_f(x)\,dx\\
		&\leq C(n,p^*,M) k^{n+\delta}2^{k\frac{n(n-{p^*})}{p^*}}\int_{B^n(0,r_{k+3})}J_f(x)\,dx\\
		&\leq C(n,p^*,M,a) k^{n+\delta}2^{k(a+\frac{n(n-{p^*})}{p^*})}.
	\end{align*}
	Proceeding with the estimation of the measure of $E(f)$ at (\ref{mofd1}):
	\[
		\begin{aligned}
			\mathcal{H}_\varepsilon^t{(E(f))}&\leq C(n,p^*,M,a)\sum_{k=M}^\infty  p_k2^{-kt}\\
			&\leq C \sum_{k=M}^\infty  k^{n+\delta}2^{-k(t-a-\frac{n(n-{p^*})}{p^*})}.
		\end{aligned}
	\]
	The sum on the right hand side converges since
	\[\frac{n(n-{p^*})}{p^*}=\frac{n}{p}~,~~t>a+\frac{n}{p}\]
	and thus we may conclude that $\mathcal{H}^t(E(f))=0$.
	
	\begin{rem}\label{rmk:on final result}
		One can show that bounded mapping $f\colon\ball\to\R^n$ of $p$-integrable distortion satisfies
		\[\int_{B^n(0,r_k)}J_f(x)\,dx\leq c(1-r_k)^{\frac{p+1}{p}(1-n)},\]
		for all $k\in\N$. On the other hand above argument shows that if
		\[\int_{B^n(0,r_k)}J_f(x)\,dx\leq c(1-r_k)^{-a},\]
		where $a\in\mathopen[0,n-1-\frac{n}{p}\mathclose)$, then $f$ has radial limits at almost every point on $\partial\ball$. So unlike in the case of quasiregular mappings, there is an interval of values of $a$ so that we don't know whether the radial limits exist a.e. on $\partial\ball$.
	\end{rem}
	
\end{proof}


\end{document}